\documentstyle[amssymb,12pt]{article}
\title{Characterizing model completeness among 
mutually algebraic structures}

\author{Michael C.\
Laskowski\thanks{Partially supported
by NSF grant DMS-0901336.}\\Department of 
Mathematics\\University of Maryland}

\newbox\smilebox
\newbox\anchorbox
\newbox\noanchorbox
\newbox\tempbox

\setbox\smilebox=\hbox{$\smile$}

\def\anchor{\hbox{\vtop{
           \hbox to \wd\smilebox{\hfil\vrule width.4pt height7pt depth1pt\hfil}
           \vskip  -11.5truept
           \hbox to \wd\smilebox{\hfil$\smile$\hfil}}}}
\setbox\anchorbox=\anchor
\def\noanchor{\hbox{\vtop{
           \hbox to \wd\anchorbox{\hfil\anchor\hfil}
           \vskip -14truept
           \hbox to \wd\anchorbox{\hfil/\hfil}}}}
\setbox\noanchorbox=\noanchor

\def\fg#1#2#3{\setbox\tempbox=\hbox{$\scriptstyle{#2}$}
\ifnum\wd\anchorbox>\wd\tempbox\dimen255=\wd\anchorbox
\else\dimen255=\wd\tempbox\fi
{#1\,\vtop{\hbox to \dimen255{\hfil\anchor\hfil}
           \vskip -6truept
           \hbox to \dimen255{\hfil$\scriptstyle{#2}$\hfil}}
           \,#3}}

\def\nfg#1#2#3{\setbox\tempbox=\hbox{$\scriptstyle{#2}$}
\ifnum\wd\noanchorbox>\wd\tempbox\dimen255=\wd\noanchorbox
\else\dimen255=\wd\tempbox\fi
{#1\,\vtop{\hbox to \dimen255{\hfil\noanchor\hfil}
           \vskip -6truept
           \hbox to \dimen255{\hfil$\scriptstyle{#2}$\hfil}}
           \,#3}}

\setbox1=\hbox{$\bot$}

\def\north#1#2{#1\,
\hbox{$\bot$\llap {\hbox to\wd1 {\hfil $/$\hfil}}}
\,#2}

\def\nao#1#2#3{#1\  \hbox{\vtop{ 
\baselineskip=4pt
\hbox{$\bot$\llap {\hbox to\wd1 {\hfil $/$\hfil}}
\hskip .05em \llap{\hbox{$^{\scriptscriptstyle{a}}$}}}\hbox{$\scriptstyle
{#2}$}}}\, #3}

\def\bp{\par{\bf Proof.}$\ \ $}

\def\includeE#1{{\lhook\kern-3.5pt\joinrel\smash{
    \mathop{\longrightarrow}\limits^{#1}}}}

\def\efor/{Example~\ref{E4}}

\def\BL/{Baldwin--Lachlan}
\def\Bu/{Buechler}
\def\Hr/{Hrushovski}
\def\lm/{locally modular}
\def\wm/{weakly minimal}
\def\nm/{non--modular}
\def\ss/{superstable}
\def\ud/{unidimensional}
\def\sm/{strongly minimal}

\def\abar{\overline{a}}

\def\ebar{\overline{e}}

\def\hbar{\overline{h}}

\def\ubar{\overline{u}}
\def\vbar{\overline{v}}
\def\wbar{\overline{w}}
\def\xbar{\overline{x}}
\def\ybar{\overline{y}}
\def\zbar{\overline{z}}

\def\acl{{\rm acl}}

\def\tp{{\rm tp}}

\def\tr/{trivial}
\def\nt/{non--trivial}
\def\st/{strong type}

\def\TV/{Tarski--Vaught}

\def\sc/{sound construction}
\def\ac/{atomic construction}

\def\fal/{functional}
\def\upl/{unique parallel lines}
\def\chp/{categorical in a higher power}

\def\conc{{\char'136}}
\def\abar{\bar{a}}

\def\ebar{\bar{e}}
\def\xbar{\bar{x}}
\def\ybar{\bar{y}}
\def\zbar{\bar{z}}
\def\phi{\varphi}

\def\A{{\cal A}}
\def\E{{\cal E}}

\def\F{{\cal F}}
\def\P{{\cal P}}

\def\MA{{\cal MA}}
\def\tp{{\rm tp}}

\def\range{{\rm range}}

\def\acl{{\rm acl}}

\def\bp{{\bf Proof.}\quad}
\def\endproof{\medskip}
\def\<{\langle}
\def\>{\rangle}
\newtheorem{Theorem}{Theorem}[section]

\newtheorem{Definition}[Theorem]{Definition}

\newtheorem{Lemma}[Theorem]{Lemma}

\begin{document}
\maketitle
\begin{abstract}
We characterize when the elementary diagram of a mutually algebraic structure has a   model complete
theory,
and give an explicit description of a set of existential formulas to which every formula is equivalent.
This characterization yields a new, more constructive proof that the elementary diagram
of any model of a
strongly minimal, trivial theory is model complete.
\end{abstract}

\section{Introduction}

In \cite{MCLma}, which borrows heavily from \cite{MCLwm},
it is shown that for any mutually algebraic structure $M$ (see Definition~\ref{mastruct}),
its elementary diagram, which we denote by $T(M)$, has a near model complete theory.
Indeed, Definition~\ref{e} describes a specific class $\E$ of existential $L(M)$-formulas, and
every $L(M)$-formula is $T(M)$-equivalent to some boolean combination of formulas from $\E$.

In earlier papers, it was shown that under stronger hypotheses on the theory
of $M$, 
 the elementary diagram $T(M)$
has a model complete theory.   Indeed, in \cite{5}, Goncharov, Harizanov, Lempp,
McCoy, and the author prove that the elementary diagram of
every model of a strongly minimal, trivial theory is model complete.
In \cite{DLR}, this result was strengthened by Dolich, Raichev, and the author
to give the same result for any model of an $\aleph_1$-categorical, trivial
theory of Morley rank 1.
In both instances,  it follows that every $L(M)$-formula is
equivalent to an existential formula, but the proofs do not give a specific 
description of a `minimal set' of 
existential formulas needed to describe all $L(M)$-formulas.

The main theorem of this short note, Theorem~\ref{char}, characterizes when the
elementary diagram of a mutually
algebraic structure $M$ has a model complete theory (as opposed to simply being near model complete).
Moreover, we display a set $\P$ of easily understood existential 
formulas\footnote{Every $\phi(\ybar)\in\P$ can be written in the form
$\exists\xbar\psi(\xbar,\ybar)$, where $\psi$ is quantifier free and
there is an integer $K$ so that
$T(M)\models\forall\ybar\exists^{<K}\xbar\psi(\xbar,\ybar)$.
Perhaps such a formula should be called an `algebraically existential'
formula?},
and show
that  $T(M)$ is model complete if and only if every $L(M)$-formula is $T(M)$-equivalent to an
element of $\P$.  Then, in the third section, we indicate that these conditions hold
for models of either of the two types of theories described above.

We conclude the Introduction by recalling the major definitions and results from \cite{MCLwm} and \cite{MCLma}.

\begin{Definition}  {\em 
When we write a tuple $\zbar$ of variable symbols, 
we assume that the elements of $\zbar$ are distinct, and 
$\range(\zbar)$ denotes
the underlying set of variable symbols.
A {\em proper partition}
$\zbar=\xbar\conc\ybar$ satisfies $\lg(\xbar),\lg(\ybar)\ge 1$, 
$\range(\xbar) \cup\range(\ybar)=\range(\zbar)$, and 
$\range(\xbar)\cap\range(\ybar)=\emptyset$.  We do not require $\xbar$ be
an initial segment of $\zbar$ 
but to simplify notation, we write it as if it were.
}
\end{Definition}

\begin{Definition}   \label{ma}
{\em Let $M$ denote any $L$-structure.  An $L(M)$-formula $\phi(\zbar)$
is {\em mutually algebraic} if there is an integer $N$
so that $M\models\forall\ybar\exists^{\le N}\xbar\phi(\xbar,\ybar)$ for every proper partition
$\xbar\conc\ybar$ of $\zbar$.
We let $\MA(M)$ denote the set of all mutually algebraic $L(M)$-formulas.  When $M$ is understood,
we simply write $\MA$.
}
\end{Definition}

The reader is cautioned that whether a formula $\phi(\zbar)$ is mutually algebraic or not
depends on the choice of free variables.  In particular, mutual algebraicity is
{\bf not} preserved under adjunction of dummy variables.  Note that
every $L(M)$-formula $\phi(z)$ with
exactly one free variable symbol is mutually algebraic.
Furthermore, note that inconsistent formulas are mutually algebraic.

The following Lemma indicates some of the closure properties of the set $\MA$.
In what follows, when we write $\phi(\xbar,\ybar)\in\MA$, we mean that 
$\phi(\zbar)\in\MA$ for any tuple $\zbar$ of distinct symbols such that 
$\range(\zbar)=\range(\xbar)\cup\range(\ybar)$, but that we are concentrating on a specific
proper partition $\zbar=\xbar\conc\ybar$ of $\phi(\zbar)$.

\begin{Lemma}  \label{closure}
Let $M$ be any structure in any language $L$.
\begin{enumerate}
\item  If $\phi(\zbar)\in\MA$, then $\phi(\sigma(\zbar))\in\MA$ for any permutation $\sigma$ of the
variable symbols;
\item  If $\phi(\xbar,\ybar)\in\MA$ and $\abar\in M^{\lg(\ybar)}$, 
then both  $\exists\ybar\phi(\xbar,\ybar)$ and $\phi(\xbar,\abar)\in\MA$;
\item  If $\phi(\zbar)\vdash\psi(\zbar)$ and $\psi(\zbar)\in\MA$, then $\phi(\zbar)\in\MA$;
\item  For $k\ge 1$, if $\{\phi_i(\zbar_i):i<k\}\subseteq \MA$, and $\bigcap_{i<k}\range(\zbar_i)$ is nonempty,
then $\psi(\wbar):=\bigwedge_{i<k} \phi_i(\zbar_i)\in\MA$, where $\range(\wbar)=\bigcup_{i<k}\range(\zbar_i)$;
\item  If $\phi(\xbar,\ybar)\in\MA$ and $r\in\omega$, then $\theta_r(\ybar):=\exists^{\ge r}\xbar\phi(\xbar,\ybar)\in \MA$.
\end{enumerate}
\end{Lemma}

\begin{Definition}  \label{mastruct}  {\em  
Given an arbitrary $L$-structure $M$, let
$\MA^*(M)$ denote the set of all $L(M)$-formulas
that are $T(M)$-equivalent to a boolean combination of
formulas from $\MA(M)$.  A structure $M$ is {\em mutually algebraic\/}
if $L(M)=\MA^*(M)$, i.e., every $L(M)$-formula is $T(M)$-equivalent to
a boolean combination of mutually algebraic formulas.
}
\end{Definition}

It is evident that the mutual algebraicity of a structure
is preserved under elementary equivalence.
The following is the main theorem of \cite{MCLma}.

\begin{Theorem}  \label{equiv}
The following are equivalent for any theory $T$:
\begin{enumerate}
\item  Every model of $T$ is a mutually algebraic structure;
\item  Every mutually algebraic expansion of every model of $T$
is a mutually algebraic structure;
\item  $Th((M,A))$ has the nfcp for every  $M\models T$ and every expansion $(M,A)$ by a unary predicate;
\item  Every complete extension of $T$ is weakly minimal and trivial.
\end{enumerate}
\end{Theorem}

Next, we recall four classes of $L(M)$-formulas that were introduced in \cite{MCLwm}.

\begin{Definition}  \label{e}
{\em 
Let $M$ be any $L$-structure.
\begin{itemize}
\item $\A=\{$all quantifier-free, mutually algebraic $L(M)$-formulas$\}$;
\item $\E=\{$all $L(M)$-formulas of the form $\exists\xbar\theta(\xbar,\ybar)$, where $\theta\in\A\}$
(we allow $\lg(\xbar)=0$ so $\A\subseteq\E$);
\item $\A^*=\{$all $L(M)$-formulas $T(M)$-equivalent to a Boolean combination of formulas from $\A\}$; and
\item $\E^*=\{$all $L(M)$-formulas $T(M)$-equivalent to a Boolean combination of formulas from $\E\}$.
\end{itemize}
}
\end{Definition}

The following Theorem is the main result of \cite{MCLwm} (noting that by Theorem~\ref{equiv},
if $M$ is mutually algebraic, then $Th(M)$ is weakly minimal and trivial).

\begin{Theorem}  \label{thm}
Let $M$ be any mutually algebraic structure.  Then:
\begin{enumerate}
\item Every quantifier-free $L(M)$-formula $\theta(\zbar)$
is in $\A^*$.
\item Every $L(M)$-formula is $T(M)$-equivalent to a Boolean combination of 
formulas from $\E$, i.e., $\E^*=L(M)$.
\end{enumerate}
\end{Theorem}

\section{A new class of existential formulas}
We begin this section with the central definitions of the current note.

\begin{Definition}  
{\em A formula $S(\wbar)$ is a {\em partial equality diagram\/}
if it is a boolean combination of formulas of the form $w=w'$ for various $w,w'\in\wbar$.

An $L(M)$-formula $\theta(\ybar,\zbar)$ is {\em preferred}
if it has the form
$$\exists \xbar(R(\xbar,\ybar)\wedge S(\xbar,\ybar,\zbar))$$
where $\xbar,\ybar,\zbar$ are disjoint tuples of variable symbols,
$\lg(\ybar)\ge 1$, $R(\xbar,\ybar)\in\A$, and $S(\xbar,\ybar,\zbar)$ is a partial equality diagram.

Let $\P$ denote the set of all $L(M)$-formulas that are
$T(M)$-equivalent to a {\bf positive} boolean combination of preferred formulas.
}
\end{Definition}

As the quantification in a preferred formula is only over the
mutually algebraic conjunct, it is easily checked that every $\phi(\ybar)\in\P$ is $T(M)$-equivalent
to an `algebraically
existential' formula in the sense of the footnote.

\begin{Lemma}  \label{mini}
Suppose that $M$ is an infinite, mutually algebraic structure,
$\xbar,\zbar,y$ are disjoint sequences of variable symbols, $\lg(y)=1$,
and $\{R_j(\xbar_j,y,\zbar_j):j\in J\}$ is a finite set of quantifier free,
mutually algebraic formulas where, for each $j$,  $\xbar_j\subseteq\xbar$, $\zbar_j\subseteq\zbar$,
and the variable $y$ occurs in $R_j$.  Then
$T(M)\models\forall\xbar\forall y\exists\zbar \bigwedge_{j\in J}\neg R_j(\xbar_j,y,\zbar_j)$. 
\end{Lemma}

\bp  Given such a set of formulas, choose $N\succeq M$ and $\abar, b$ from $N$.  
We will produce a tuple $\ebar$ from $N$
so  that $N\models\neg R(\abar_j,b,\ebar_j)$ for each $j\in J$.  
Say $\zbar=(z_0,\dots,z_{k-1})$.  For each $\ell<k$, let  
$J_\ell=\{j\in J:z_\ell$ occurs in $\zbar_j\}$ and
let
$$B_\ell:=\{c\in N:
 N\models\exists\zbar_j[R_j(\abar_j,b,\zbar_j)\wedge z_\ell=c]\ 
\hbox{for some $j\in J_\ell$}
\}$$
As each $R_j$ is mutually algebraic and $b$ is fixed, it follows that each of the sets $B_\ell$
is finite.  Since $N$ is infinite, we can choose $\ebar=(e_0,\dots,e_{k-1})$ so that $e_\ell\not\in B_\ell$
for each $\ell<k$.  
It is easily checked that $\ebar$ is as desired.
\endproof

\begin{Lemma}  \label{messy}
Let $M$ be an infinite, mutually algebraic structure.
Say $$\psi(\xbar,y):=\bigwedge_{i\in I} R_i(\xbar_i,y)\wedge\bigwedge_{j\in J} \neg R_j(\xbar_j,y)$$
where $I$ and $J$ are finite, each $R_i,R_j$ is quantifier free and mutually algebraic, each $\xbar_i$
and $\xbar_j$
is a subsequence of $\xbar$, $\lg(y)=1$, and $y$ occurs in each $R_i$, $R_j$.  
Then $\exists\xbar \psi(\xbar,y)\in\P$.
\end{Lemma}

\bp  First, if $I=\emptyset$, then by Lemma~\ref{mini}, $T(M)\models\forall y\exists\xbar\psi(\xbar,y)$,
hence $\exists\xbar\psi(\xbar,y)$ is true for every $y$.  In this case,  
$\exists\xbar\psi$ is equivalent to
$y=y$, which is in $\A$, and hence in $\P$.

Next, assume that $I\neq\emptyset$.  Let $\xbar'$ be the smallest subsequence of $\xbar$ for which
every $\xbar_i$ is a subseqence of $\xbar'$.  Let $\zbar=\xbar\setminus\xbar'$, let $K=\{j\in J:\xbar_j\subseteq\xbar'\}$
and let $J^*=J\setminus K$.
As $I$ is non-empty, it follows from Lemma~\ref{closure}(3) and (4) that the formula
$$\theta(\xbar',y):=\bigwedge_{i\in I} R_i(\xbar_i,y)\wedge\bigwedge_{j\in K} \neg R_j(\xbar_j,y)$$
is mutually algebraic (and it is visibly quantifier free).
But, by Lemma~\ref{mini}, it follows that $\exists\xbar\psi(\xbar,y)$ is $T(M)$-equivalent to $\exists\xbar'\theta(\xbar',y)$,
so $\exists\xbar\psi(\xbar,y)\in\P$.
\endproof

\begin{Theorem}  \label{char}
The following are equivalent for every mutually algebraic structure $M$:
\begin{enumerate}
\item  $\exists^{=r}\xbar R(\xbar,y)\in\P$ for all $R(\xbar,y)\in\A$ with $\lg(y)=1$ and all $r\in\omega$;
\item  $\exists^{=r}\xbar R(\xbar,y)\in\P$ for all $R(\xbar,y)\in\A$ with $\lg(y)\ge 1$ and all $r\in\omega$;
\item  $\P$ is closed under negation;
\item  $\P=L(M)$;
\item  $T(M)$ is model complete.
\end{enumerate}
\end{Theorem}

\bp    First, note that if the universe of $M$ is finite, then all five conditions hold trivially.
Thus, we assume throughout that $M$ is infinite.

$(1)\Rightarrow(2)$:  Assume that (1) holds.  
Choose any $R(\xbar,\ybar)\in \A$ and any integer $r$.  Choose any variable symbol
$y^*\in \ybar$ and let $\ybar'$ satisfy $\ybar'\conc y^*=\ybar$.  Choose an integer $N$ so that
$R(\xbar\ybar',y^*)$ has fewer than $N$ solutions.  
For each $m<N$, let
$$S_m(\ubar_0\vbar_0\dots\ubar_{m-1}\vbar_{m-1},\ybar'):=\bigwedge_{i\neq j} \ubar_i\vbar_i\neq\ubar_j\vbar_j\wedge
\bigvee_{Q\in{{m\choose r}}} \left(\bigwedge_{i\in Q}\vbar_i=\ybar'\wedge\bigwedge_{i\not\in Q}\vbar_i\neq\ybar'\right)$$
and let $$\theta_m(\ybar):=
\exists\ubar_0\vbar_0\dots\exists\ubar_{m-1}\vbar_{m-1}\left(\bigwedge_{i<m} R(\ubar_i\vbar_i,y^*)\wedge S_m(\ubar_0\vbar_0\dots
\ubar_{m-1}\vbar_{m-1}\ybar')\right).$$
Using the closure properties in Lemma~\ref{closure}, $\theta_m(\ybar)$ is a preferred formula.
Let $\wbar$ be new variables satisfying $\lg(\wbar)=\lg(\ubar)+\lg(\vbar)$ and let $\delta(\ybar)$ be
$$\bigvee_{m<N}\left(\exists^{=m}\wbar R(\wbar,y^*)\wedge
\theta_m(\ybar)\right)$$
It is easily checked that $\delta(\ybar)$ is $T(M)$-equivalent to 
$\exists^{=r}\xbar R(\xbar,\ybar)$ and, using (1), $\delta(\ybar)\in\P$.

$(2)\Rightarrow(3)$:
In order to show that $\P$ is closed under negation, by DeMorgan's laws
it suffices to show that the negation of every preferred formula
is in $\P$.  
So fix a preferred formula $\theta(\ybar,\zbar):=\exists\xbar(R(\xbar,\ybar)\wedge S(\xbar,\ybar,\zbar))$,
where $R(\xbar,\ybar)\in\A$, $\lg(\ybar)\ge 1$, and $S(\xbar,\ybar,\zbar)$ is a partial equality diagram.
Choose $N$ so that $T(M)$ implies that $\exists^{<N}\xbar R(\xbar,\ybar)$.
It is easily checked that $\neg\theta(\ybar,\zbar)$ is $T(M)$-equivalent to
$$\bigvee_{m<N} \left(\exists^{=m}\xbar R(\xbar,\ybar)\wedge\psi_m(\ybar,\zbar)\right)$$
where $$\psi_m(\ybar,\zbar):=\exists\xbar_0\dots\xbar_{m-1}\left(\bigwedge_{i<m} R(\xbar_i,\ybar)\wedge
\bigwedge_{i\neq j} \ybar_i\neq\ybar_j\wedge\bigwedge_{i<m}\neg S(\xbar_i,\ybar,\zbar)\right)$$
Thus, $\neg\theta(\ybar,\zbar)\in\P$ by (2).

$(3)\Rightarrow(4)$:  As $\P$ is closed under positive boolean combinations by definition, it follows
immediately from (3) that $\P$ is closed under {\bf all} boolean combinations.
However, $\E\subseteq \P$ trivially, so $\E^*$, the closure of $\E$ under boolean combinations, is
also a subset of $\P$.  But, as $M$ is mutually algebraic, $\E^*=L(M)$ by Theorem~\ref{thm}(2).  
Thus $\P=L(M)$.

$(4)\Rightarrow(5)$:  Visibly, every preferred formula is an existential $L(M)$-formula, and the set of
existential $L(M)$-formulas is closed under positive
boolean combinations.  Thus, (4) implies that every $L(M)$-formula
is $T(M)$-equivalent to an existential formula, which is equivalent to model completeness (see e.g., \cite{CK}).

$(5)\Rightarrow(1)$:  Assume that $T(M)$ is model complete.  We argue that every $L(M)$-formula
$\phi(y)$ with $\lg(y)=1$ is in $\P$.  Fix such a formula $\phi(y)$.

\medskip
\par\noindent{\bf Claim:}  For any $N\succeq M$ and any $b\in N$ such that $N\models\phi(b)$,
there is $\delta(y)\in\P$ such that $N\models\delta(b)\wedge\forall y(\delta(y)\rightarrow\phi(y))$.

\medskip

\bp  Fix such an $N$ and $b$.
As $T(M)$ is model complete, this  implies that 
$$T(M)\cup \Delta_{M^*}\models\phi(b)$$
where $\Delta_{M^*}$ is the atomic diagram of $M^*$.  Thus, by compactness, there is
a quantifier-free  $\theta(\ebar,b)\in\Delta_{M^*}$ such that $T(M)\cup\{\theta(\ebar,b)\}\models\phi(b)$.
Without loss, we may assume that $\ebar$ is disjoint from $M\cup\{b\}$, so it follows
that
$T(M)\models\forall y(\exists\xbar\theta(\xbar,y)\rightarrow\phi(y))$.

By Theorem~\ref{thm}(1), $\theta(\xbar,y)\in\A^*$.  Thus, by considering the Disjuntive Normal Form,
we can write $\theta(\xbar,y)$ as $\bigvee\bigwedge R_{ij}(\zbar_{ij})$, where each 
each $\zbar_{ij}$ is contained in $\xbar\cup\{y\}$ and $R_{ij}(\zbar_{ij})$
quantifier-free and 
is either mutually algebraic or the is negation of a mutually algebraic formula. 

Thus, one of the disjuncts $\psi(\xbar',y)$ of $\theta(\xbar,y)$ satisfies
$N\models\exists\xbar'\psi(\xbar',b)$, 
$$T(M)\models\forall y(\exists\xbar'\psi(\xbar',y)\rightarrow\phi(y))$$
and $\xbar'\subseteq\xbar$.
Now $\psi(\xbar',y)$ has the form
$$\bigwedge R_i(\xbar_i,y)\wedge\bigwedge \neg R_j(\xbar_j,y)$$
where each $R_i$ and $R_j$ is quantifier-free and mutually algebraic and each
$\xbar_i,\xbar_j\subseteq\xbar'$.
We may additionally assume that the variable symbol $y$ appears in each $R_i$ and $R_j$.
As $M$ is infinite, Lemma~\ref{messy} applies, 
and the formula $\delta(y):=\exists\xbar'\psi(\xbar',y)\in\P$
is as required.
\endproof

To finish the proof of $(5)\Rightarrow(1)$,
let $$\Gamma:=\{\delta(y)\in\P:T(M)\models\forall y(\delta(y)\rightarrow\phi(y))\}$$
It follows immediately from the Claim and compactness
that the formula $\phi(y)$ is $T(M)$-equivalent to
a finite disjunction 
$\bigvee_i\delta_i(y)$ of elements $\delta_i\in\Gamma$.  As $\P$ is closed under 
$T(M)$-equivalence and finite disjunctions
we conclude that $\phi(y)\in\P$.
\endproof

\section{New proofs of model completeness}  

We close by giving new proofs of the model completeness results first proved in
\cite{5} and \cite{DLR}.  The first theorem clearly follows from the second, but
we give a separate proof as it follows so easily from our main result.

\begin{Theorem}  \label{sm}  
If $T$ is strongly minimal and trivial, then $T(M)$ is model complete and
$L(M)=\P$ for every
model $M$ of $T$.
\end{Theorem}

\bp  Fix a model $M$ of $T$.   With our eye on Clause (1) of Theorem~\ref{char},
choose an $L(M)$-formula $\phi(y)$ with $\lg(y)=1$.
 By strong minimality, the solution set 
$\phi(N)$ in any $N\succeq M$ is either finite
or cofinite, with the `exceptional set' contained in $M$.  That is,
there is some finite set $Q\subseteq M$ such that, letting 
$\theta(y):=\bigvee_{m\in Q} y=m$, $\phi(y)$
is $T(M)$-equivalent to either $\theta(y)$ or $\neg\theta(y)$.  As any 
quantifier free $L(M)$-formula in a single free variable is in $\A$ and hence
in $\P$, both
$\theta,\neg\theta\in\P$.  Applying this argument  to any instance
of $\exists^{=r}\zbar R(\zbar,y)$, we conclude that both $T(M)$ is model
complete and $L(M)=\P$ by Theorem~\ref{char}.
\endproof

\begin{Theorem}  \label{rank1}
Suppose $T$ is $\aleph_1$-categorical, trivial, and of Morley rank 1.
Then for every $M\models T$, the elementary diagram is model complete.
Furthermore, $L(M)=\P$.
\end{Theorem}

\bp  Again, we employ Theorem~\ref{char}, but here we need to
focus on a particular instance of Clause~(1).  
So fix a formula
$R(\zbar,y)\in \A$ and  an integer $r$.  As $\exists^{\ge r}\zbar R(\zbar,y)\in\P$,
to establish Clause~(1) it suffices to prove that $\exists^{\le r}\zbar R(\zbar,y)\in\P$.  

Toward this end, our assumptions on $T$ imply that there are finitely
many non-algebraic 1-types over $M$.  Indeed, if $S_{na}:=\{p_i:i<d\}$
denotes this set of non-algebraic 1-types, then 
$d$ is the Morley degree of $T$.  As well, 
the $\aleph_1$-categoricity of $T$ implies each of these types are non-orthogonal.
As $T$ is trivial, this further implies that $\nao {p_i} {M} {p_j}$ for all
$p_i,p_j\in S_{na}$.  As forking of a 1-type implies algebraicity,
this implies that for any $N\succeq M$ and any
$a\in p_i(N)$, there is $b\in p_j(N)$
such that $b\in\acl(M\cup\{a\})$ (and hence $a\in\acl(M\cup\{b\})$).
As $\E^*=L(M)$, it is easy to verify that for all pairs $p_i,p_j\in S_{na}$,
there is a mutually algebraic, quantifier free formula $\theta_{ij}(x,y,\zbar)$
such that for any $a\in p_i(N)$, there is $b\in p_j(N)$ such that
$N\models\exists\zbar\theta_{ij}(a,b,\zbar)$.
Fix a finite set $\F\subseteq\A$ consisting of one such $\theta_{ij}$ for
each pair $p_i,p_j\in S_{na}$ (if $i=j$ we can take $\theta_{ij}$ to be
the mutually algebraic
formula $x=y$).

Now, fix an elementary extension $N\succeq M$.
For every $b\in N$, let
\begin{itemize}
\item $\theta^*(b,N)=\{c\in N:N\models\exists\wbar\theta(b,c,\wbar)\}$
i.e., $\theta^*(b,N)$ is the {\bf set} of elements that are part of a tuple realizing
$\theta(b,N)$;
\item  $\F(b)=\bigcup_{\theta\in \F} \theta^*(b,N)$; and
\item $\F_r(b)=\{c\in \F(b):N\models\exists^{\ge r+1}\zbar R(\zbar,c)\}$.
\end{itemize}

Clearly, for $b\in N\setminus M$, $\F_r(b)\subseteq\F(b)$ and
$|\F(b)|\le\sum_{\theta\in \F} N_\theta\cdot\lg(\zbar)$, 
where $N_\theta$ is an integer
such that $T(M)\models\forall y\exists^{<N_\theta}\zbar\theta(y,\zbar)$.

Thus, there is a finite exceptional set $Q\subseteq M$ and an integer $\ell^*$
such that 
\begin{enumerate}
\item $T(M)\models\forall y (y\not\in Q\rightarrow |\F_r(y)|\le \ell^*)$ and
\item  For some $b\in N\setminus M$, $|\F_r(b)|=\ell^*$.
\end{enumerate}

Also, it is clear that the size $|\F_r(b)|$ depends only on $\tp(b/M)$, i.e.,
if $\tp(b/M)=\tp(b'/M)$, then $|\F_r(b)|=|\F_r(b')|$.
Fix any non-algebraic 1-type $p^*(y)\in S_{na}$ such that
$|\F_r(b)|=\ell^*$ for some (every) realization $b$ of $p^*$.

Let $\delta(x)$ express
\begin{quotation}
``There is some $\theta(x,y,\zbar)\in\F$ such that
$\exists y\exists\zbar \bigg(\theta(x,y,\zbar)\wedge y\not\in Q$
and
there are distinct elements $\{w_i:i<\ell^*\}$ witnessing that $|\F_r(y)|\ge\ell^*$
{\bf and} $x\neq w_i$ for all $i<\ell^*\bigg)$.''
\end{quotation}
It is routine to check that the formula $\delta(x)\in\P$.
It suffices to prove the following:

\medskip\par\noindent
{\bf Claim:}  $T(M)\models\forall x[\exists^{\le r}\zbar R(\zbar,x)\leftrightarrow
\delta(x)]$.

\medskip

\bp  Fix any $N\succeq M$ and $a\in N$. First,
 suppose $N\models\exists^{\le r} \zbar R(\zbar,a)$.  Choose $\theta(x,y,\zbar)\in\F$
such that there is $b\in p^*(N)$ with $N\models\exists\zbar\theta(a,b,\zbar)$.
By our choice of $p^*$ we have $|\F_r(b)|=\ell^*$, so choose an enumeration
$\{c_i:i<\ell^*\}$ of $\F_r(b)$.  Since the definition of $\F_r(b)$ implies
that $N\models\exists^{\ge r+1}\zbar R(\zbar,c_i)$ for each $i$, it follows
that $a\neq c_i$ for each $i$.  Thus, $N\models\delta(a)$.

Conversely, suppose that $N\models\exists^{\ge r+1}\zbar R(\zbar,a)$.
Choose any $b\in N\setminus Q$ such that $N\models\exists\zbar\theta(a,b,\zbar)$
for some $\theta\in\F$ and $|\F_r(b)|\ge\ell^*$.
Choose any set of $\ell^*$ distinct elements
$\{c_i:i<\ell^*\}\subseteq\F_r(b)$.  But now, as $b\not\in Q$,
$$N\models |\F_r(b)|\le\ell^*$$
This, combined with the fact that our assumption on $a$ and $\theta$ imply
that $a\in\F_r(b)$, guarantees that $a=c_i$ for some $i$.
That is, $N\models\neg\delta(a)$, completing the proof of the Claim.
\endproof

As we have shown that $\exists^{\le r}\zbar R(\zbar,x)\in\P$, it follows
from Theorem~\ref{char} that both $T(M)$ is model complete and
$L(M)=\P$.
\endproof

\end{document}